\documentclass[12pt]{amsart}

\usepackage{amsmath,amssymb,amscd,amsthm}
\usepackage{hyperref}
\usepackage[left=3cm,right=3cm,top=3cm,bottom=3cm,includeheadfoot]{geometry}

\makeatletter
\@namedef{subjclassname@2010}{%
  \textup{2010} Mathematics Subject Classification}
\makeatother

\newtheorem{theorem}{Theorem}

\theoremstyle{definition}

\DeclareMathOperator\euler{Euler}

\begin{document}

\baselineskip=16pt

\title[The algebraic degree in semidefinite programming]{A formula for the algebraic degree\\ in semidefinite programming}

\author{Dang Tuan Hiep}

\address{National Center for Theoretical Sciences\\ No. 1 Sec. 4 Roosevelt Rd., National Taiwan University\\ Taipei, 106, Taiwan}

\email{hdang@ncts.ntu.edu.tw}

\begin{abstract}
In this paper we use the Bott residue formula in equivariant cohomology to show a formula for the algebraic degree in semidefinite programming.
\end{abstract}



\date{\today}

\maketitle

\section{Introduction}

Consider the semidefinite programming (SDP) problem in the form 
\begin{equation}\label{SDP}
  \text{maximize trace}(B \cdot Y) \text{ subject to } Y \in \mathcal U \text{ and } Y \succeq 0,
\end{equation}
where $B$ is a real symmetric $n\times n$-matrix, $\mathcal U$ is a $m$-dimensional affine subspace in the $\binom{n+1}{2}$-dimensional space of real $n \times n$-symmetric matrices, and $Y \succeq 0$ means that $Y$ is positive semidefinite. We know that the coordinates of the optimal solution are the roots of some univariate polynomials. If the data are generic, then the degree of these polynomials depends only on the rank $r$ of the optimal solution. This is what we call the {\it algebraic degree} $\delta(m,n,r)$ of the semidefinite programming \eqref{SDP}. 

Let $m,n$, and $r$ be positive integers satisfying the following conditions
\begin{equation}
  \binom{n-r+1}{2} \leq m \text{ and } \binom{r+1}{2} \leq \binom{n+1}{2} - m.
\end{equation}
Nie, Ranestad, and Sturmfels \cite[Theorem 13]{NRS} showed that the number $\delta(m,n,r)$ is equal to the degree of the dual of a determinantal variety. In the language of intersection theory, von Bothmer and Ranestad \cite[Proposition 4.1]{BR} showed that
\begin{equation}\label{schubert}
  \delta(m,n,r) = \int_{G(r,n)}s_{m-\binom{n-r+1}{2}}(S^2Q)s_{\binom{n+1}{2}-m-\binom{r+1}{2}}(S^2U^*),
\end{equation}
where $U$ and $Q$ are respectively the universal sub-bundle and quotient bundle on the Grassmannian $G(r,n)$, $S^2Q$ and $S^2U^*$ are respectively the second symmetric power of $Q$ and the dual of $U$, and $s_i(E)$ is the $i$-th Segre class of the dual of the vector bundle $E$. Note that $\int_X\alpha$ denotes the degree of the cycle class $\alpha$ on $X$ defined in \cite[Definition 1.4]{F}. Furthermore, a general formula for $\delta(m,n,r)$ was given in \cite[Theorem 1.1]{BR}. In this paper we show another formula for $\delta(m,n,r)$. Our method is completely different from the previous one. We use the Bott residue formula in equivariant cohomology.
\subsection{Main result}
Consider the polynomial ring $\mathbb Q[\lambda_1,\ldots,\lambda_n]$ in $n$ variables $\lambda_1,\ldots,\lambda_n$. For each subset $I \subseteq \{1,\ldots,n\}$ and positive integer $i$, we define $c_{i,I}$ to be the $i$-th elementary symmetric polynomial in $\binom{|I|+1}{2}$ variables which are the elements of the set
$$\left\{\sum_{i\in I}a_i\lambda_i \mid a_i \in \{0,1,2\}, \sum_{i\in I}a_i = 2\right\},$$
$$A_{i,I} = \det\left(\begin{matrix}
c_{1,I} & c_{2,I} & c_{3,I} & \cdots & c_{i,I}\\
1 & c_{1,I} & c_{2,I} & \cdots & c_{i-1,I}\\
0 & 1 & c_{1,I} & \cdots & c_{i-2,I}\\
\vdots & \vdots & \vdots & \ddots & \vdots \\
0 & 0 & 0 & \cdots & c_{1,I}
\end{matrix}\right),$$
$$T_I = \prod_{i\in I} \prod_{j\not\in I}(\lambda_j - \lambda_i).$$
We also define $A_{0,I} = 1$ for all $I$. Note that $c_{i,I} = 0$ whenever $i > \binom{|I|+1}{2}$.
\begin{theorem}\label{main}
The algebraic degree
$$\delta(m,n,r) = (-1)^{l}\sum_I\frac{A_{k,I^c}A_{l,I}}{T_I},$$
where the sum runs over all subsets $I$ consisting of $r$ elements of $\{1,\ldots,n\}$, $I^c$ is the complement of $I$ in $\{1,\ldots,n\}$, and $k,l$ stand for $m-\binom{n-r+1}{2}, \binom{n+1}{2}-m-\binom{r+1}{2}$ respectively.
\end{theorem}
The right-hand side of the formula in Theorem \ref{main} is a rational polynomial function, and the theorem claims in other words that it is in fact a constant function, moreover it is an integer. Namely, for any numbers $\lambda_i$ such that $\lambda_i \neq \lambda_j$ for $i \neq j$, the right-hand side of the formula becomes the same integer.



\section{Proof of Theorem \ref{main}}
We use the Bott residue formula in equivariant cohomology to prove the formula in Theorem \ref{main}. For more details on the Bott residue formula, we refer to \cite[Proposition 9.1.5]{CK} for a topological version and \cite[Theorem 3]{EG} for an algebraic version. Our primary interest is as follows. Consider the diagonal action of $T = (\mathbb C^*)^n$ on $\mathbb C^n$ given in coordinates by
$$(t_1,\ldots,t_n) \cdot (x_1, \ldots , x_n) = (t_1x_1, \ldots , t_nx_n).$$
This induces an action of $T$ on the Grassmannian $G(r,n)$ with $\binom{n}{r}$ isolated fixed points $L_I$ corresponding to $\binom{n}{r}$ coordinate $r$-planes in $\mathbb C^n$. Each fixed point $L_I$ is indexed by a subset $I$ of length $r$ of the set $\{1,\ldots,n\}$.

By (\ref{schubert}) and the Bott residue formula \cite[Proposition 9.1.5]{CK}, we obtain
$$\delta(m,n,r) = \sum_{I} \frac{s^T_k(S^2Q|_{L_I})s^T_l(S^2U^*|_{L_I})}{\euler_T(N_{L_I})},$$
where the sum runs over all subsets $I$ consisting of $r$ elements of the set $\{1,\ldots,n\}$. For each $I$, $s^T_k(S^2Q|_{L_I})$, $s^T_l(S^2U^*|_{L_I})$, and $\euler_T(N_{L_I})$ are evaluated as follows.

Let $U$ and $Q$ be the universal sub-bundle and quotient bundle on $G(r,n)$ respectively. At each $L_I$, the torus action on the fibers $U|_{L_I}$ and $Q|_{L_I}$ have characters $\rho_i$ for $i\in I$ and $\rho_j$ for $j\not\in I$ respectively. Since the tangent bundle on the Grassmannian is isomorphic to $U^* \otimes Q$, the characters of the torus action on the tangent space at $L_I$ are
$$\{\rho_j - \rho_i \mid i \in I, j \not \in I \}.$$
The normal bundle $N_{L_I}$ of $L_I$ in $G(r,n)$ is just the tangent space of $G(r,n)$ at $L_I$. Hence
\begin{align*}
\euler_T(N_{L_I}) & = \prod_{i\in I} \prod_{j\not\in I} (\lambda_j-\lambda_i)\\
& = T_I.
\end{align*}
Note that the $\rho_i$ and $\lambda_i$ are defined in \cite[Subsection 9.1.1]{CK}.

Since the characters of the torus action on $U^*|_{L_I}$ are $-\rho_i$ for $i\in I$, the torus action on $S^2U^*|_{L_I}$ has the characters 
$$\left\{-\sum_{i\in I} a_i\rho_i \mid a_i\in \{0,1,2\}, \sum_{i\in I} a_i = 2 \right\}.$$
Thus $c_i^T(S^2U^*|_{L_I})$ is the $i$-th elementary symmetric polynomial in $\binom{r+1}{2}$ variables which are the elements of the set
$$\left\{-\sum_{i\in I} a_i\lambda_i \mid a_i\in \{0,1,2\}, \sum_{i\in I} a_i = 2 \right\}.$$
This implies that $c_i^T(S^2U^*|_{L_I}) = (-1)^ic_{i,I}$ for all $i$. Hence
\begin{align*}
s_l^T(S^2U^*|_{L_I}) & =  \det\left(\begin{matrix}
c_1^T(S^2U^*|_{L_I}) & c_2^T(S^2U^*|_{L_I}) & c_3^T(S^2U^*|_{L_I}) & \cdots & c_l^T(S^2U^*|_{L_I})\\
1 & c_1^T(S^2U^*|_{L_I}) & c_2^T(S^2U^*|_{L_I}) & \cdots & c_{l-1}^T(S^2U^*|_{L_I})\\
0 & 1 & c_1^T(S^2U^*|_{L_I}) & \cdots & c_{l-2}^T(S^2U^*|_{L_I})\\
\vdots & \vdots & \vdots & \ddots & \vdots \\
0 & 0 & 0 & \cdots & c_1^T(S^2U^*|_{L_I})
\end{matrix}\right)\\
& = (-1)^l\det\left(\begin{matrix}
c_{1,I} & c_{2,I} & c_{3,I} & \cdots & c_{l,I}\\
1 & c_{1,I} & c_{2,I} & \cdots & c_{l-1,I}\\
0 & 1 & c_{1,I} & \cdots & c_{l-2,I}\\
\vdots & \vdots & \vdots & \ddots & \vdots \\
0 & 0 & 0 & \cdots & c_{1,I}
\end{matrix}\right)\\
& = (-1)^lA_{l,I}.
\end{align*}
Similarly, since the characters of the torus action on $Q|_{L_I}$ are $\rho_i$ for $i \in I^c$, the torus action on $S^2Q|_{L_I}$ has the characters
$$\left\{\sum_{i\in I^c} a_i\rho_i \mid a_i\in \{0,1,2\}, \sum_{i\in I^c} a_i = 2 \right\}.$$
Thus $c_i^T(S^2Q|_{L_I})$ is the $i$-th elementary symmetric polynomial in $\binom{n-r+1}{2}$ variables which are the elements of the set
$$\left\{\sum_{i\in I^c} a_i\lambda_i \mid a_i\in \{0,1,2\}, \sum_{i\in I^c} a_i = 2 \right\}.$$
This implies that $c_i^T(S^2Q|_{L_I}) = c_{i,I^c}$ for all $i$. Therefore, we have
\begin{align*}
s_k^T(S^2Q|_{L_I}) & =  \det\left(\begin{matrix}
c_1^T(S^2Q|_{L_I}) & c_2^T(S^2Q|_{L_I}) & c_3^T(S^2Q|_{L_I}) & \cdots & c_k^T(S^2Q|_{L_I})\\
1 & c_1^T(S^2Q|_{L_I}) & c_2^T(S^2Q|_{L_I}) & \cdots & c_{k-1}^T(S^2Q|_{L_I})\\
0 & 1 & c_1^T(S^2Q|_{L_I}) & \cdots & c_{k-2}^T(S^2Q|_{L_I})\\
\vdots & \vdots & \vdots & \ddots & \vdots \\
0 & 0 & 0 & \cdots & c_1^T(S^2Q|_{L_I})
\end{matrix}\right)\\
& = \det\left(\begin{matrix}
c_{1,I^c} & c_{2,I^c} & c_{3,I^c} & \cdots & c_{k,I^c}\\
1 & c_{1,I^c} & c_{2,I^c} & \cdots & c_{k-1,I^c}\\
0 & 1 & c_{1,I^c} & \cdots & c_{k-2,I^c}\\
\vdots & \vdots & \vdots & \ddots & \vdots \\
0 & 0 & 0 & \cdots & c_{1,I^c}
\end{matrix}\right)\\
& = A_{k,I^c}.
\end{align*}
In summary, we obtain the desired formula.

In the same way, we also obtain a formula \cite[Theorem 1.1]{H1} for the degree of Fano schemes of linear subspaces on hypersurfaces. An implementation of the formula in Theorem \ref{main} has been done in \texttt{Schubert3} \cite{H2}, which is a \textsc{Sage} \cite{S} package for computations in intersection theory and enumerative geometry. 

\subsection*{Acknowledgements}
Part of this work was done while the author was visiting the University of Kaiserslautern and Oberwolfach Research Institute for Mathematics. This work is finished during the author's postdoctoral fellowship at the National Center for Theoretical Sciences. He thanks all for the financial support and hospitality.

\end{document}